# SPATIAL EXTREMES: MODELS FOR THE STATIONARY CASE

## By Laurens de Haan[1] and Teresa T. Pereira[2]

### *Erasmus University and University of Lisbon*

The aim of this paper is to provide models for spatial extremes in the case of stationarity. The spatial dependence at extreme levels of a stationary process is modeled using an extension of the theory of max-stable processes of de Haan and Pickands [*Probab. Theory Related Fields* **72** (1986) 477–492]. We propose three one-dimensional and three two-dimensional models. These models depend on just one parameter or a few parameters that measure the strength of tail dependence as a function of the distance between locations. We also propose two estimators for this parameter and prove consistency under domain of attraction conditions and asymptotic normality under appropriate extra conditions.

**1. Introduction.** The paper develops a framework as well as concrete models for statistics of spatial extremes which are sufficiently simple to be used in applications. Only the case of stationary processes is considered and the dependence structure will be represented by one parameter or a few parameters. Instead of developing more complicated models we aim at developing several simple models with somewhat different features.

For simplicity of exposition and in order to stay close to the existing literature we shall start discussing processes which are defined on $\mathbb{R}$ rather than $\mathbb{R}^2$. After that we discuss stationary processes in $\mathbb{R}^2$.

The setting is as follows. Consider independent replications of a stochastic process with continuous sample paths

$$\{X_n(t)\}_{t \in \mathbb{R}},$$

Received November 2002; revised March 2005.

[1]Supported in part by Gulbenkian Foundation.

[2]Supported by FCT/POCTI/35163/MAT/2000/FEDER (Project MEDE—*Statistical Modeling of Spatial Data*).

$n = 1, 2, \ldots$. Suppose that the process is in the domain of attraction of a max-stable process, that is, there are sequences of continuous functions $a_n > 0$ and $b_n$ such that as $n \to \infty$

$$(1.1) \qquad \left\{ \frac{\max_{1 \leq i \leq n} X_i(t) - b_n(t)}{a_n(t)} \right\}_{t \in \mathbb{R}} \xrightarrow{W} \{\tilde{Z}(t)\}_{t \in \mathbb{R}}$$

in $C$-space. Necessary and sufficient conditions have been given by de Haan and Lin [4]. The limit process $\{\tilde{Z}(t)\}$ is a max-stable process. Without loss of generality we can assume that the marginal distributions of $\tilde{Z}$ can be written as

$$\exp\{-(1 + \gamma(t)x)^{-1/\gamma(t)}\}$$

for all $x$ with $1 + \gamma(t)x > 0$ where the function $\gamma$ is continuous. For the time being we shall discuss the standardized process, called simple max-stable,

$$\{Z(t)\}_{t \in \mathbb{R}} := \{(1 + \gamma(t)\tilde{Z}(t))_+^{1/\gamma(t)}\}_{t \in \mathbb{R}},$$

whose marginal distribution functions are all Fréchet: $\exp(-1/x)$, $x > 0$.

$\{Z(t)\}$ is a simple max-stable process. We assume that $\{Z(t)\}$ is a stationary process. The theory of de Haan and Pickands [5] applies. According to Theorem 6.1 of that paper the process is determined by a nonnegative $L_1$ function and a group of linear $L_{1+}$ isometries. However, since we aim at manageable models, we shall restrict ourselves to the subclass of stationary max-stable processes which is discussed on pages 490–491 of [5], the one of moving maximum processes. The process is defined as follows.

Let $\phi$ be a unimodal continuous probability density on the real line and let $\{X_j, Y_j\}_{j \geq 1}$ be the points of a homogeneous Poisson process on $\mathbb{R} \times \mathbb{R}_+$. The process is defined as a functional of the Poisson process as follows:

$$Z(t) := \max_{j \geq 1} \frac{\phi(X_j - t)}{Y_j} \qquad \text{for } t \in \mathbb{R}.$$

It is easy to check directly that this process is stationary and simple max-stable. The almost sure continuity of the process follows from [1]. We think of $t$ as a space parameter, not a time parameter. Note that for $t_1, t_2, \ldots, t_d \in \mathbb{R}$ and $x_1, x_2, \ldots, x_d > 0$ (cf. [5]),

$$P\{Z(t_1) \leq x_1, \ldots, Z(t_d) \leq x_d\} = \exp\left\{ -\int_{-\infty}^{+\infty} \max_{1 \leq i \leq d} \frac{\phi(s - t_i)}{x_i} \, ds \right\}.$$

However, this is not yet sufficiently simple for applications. We shall consider three specific examples depending on just one parameter. For $\phi$ we choose the normal density

$$(1.2) \qquad \frac{\beta}{\sqrt{2\pi}} \exp\left\{ -\frac{\beta^2 x^2}{2} \right\},$$

the double exponential density

$$(1.3) \qquad\qquad \frac{\beta}{2} \exp\{-\beta|x|\}$$

and the $t$-density

$$(1.4) \qquad \frac{\beta\Gamma((\nu+1)/2)}{\sqrt{\pi\nu}\Gamma(\nu/2)} \left\{1 + \frac{\beta^2 x^2}{\nu}\right\}^{-(\nu+1)/2}, \qquad \text{with } \nu \text{ a positive integer,}$$

where $\beta$ is a positive constant. The constant $\beta$ measures the strength of tail dependence. We shall see that in all cases small values of $\beta$ point at strong dependence and large values of $\beta$ point at weak dependence. However, the dependence does not decrease at the same rate in all cases.

As it should be for spatial models, the tail dependence between $Z(0)$ and $Z(t)$ decreases monotonically and continuously with $|t|$. In particular, when $|t| \to \infty$, the random variables $Z(0)$ and $Z(t)$ become independent.

The same happens for fixed $t$ when varying $\beta$: as $\beta \downarrow 0$ the process becomes a.s. constant and as $\beta \to \infty$ $Z(0)$ and $Z(t)$ become independent (cf. propositions below). The dependence decreases monotonically as $\beta$ increases.

Next we extend the definition to processes on $\mathbb{R}^2$, that is, to random fields. It is readily seen that the theory of de Haan and Pickands [5] remains valid if the underlying Poisson process is based on $\mathbb{R}^2 \times \mathbb{R}_+$ rather than $\mathbb{R} \times \mathbb{R}_+$. So we consider a unimodal (i.e., nonincreasing in each direction starting from the mode) continuous probability density $\phi$ on $\mathbb{R}^2$. Let $\{X_j, W_j, Y_j\}_{j \geq 1}$ be the points of a homogeneous Poisson point process on $\mathbb{R}^2 \times \mathbb{R}_+$. The process defined by

$$Z(t_1, t_2) := \max_{j \geq 1} \frac{\phi(X_j - t_1, W_j - t_2)}{Y_j} \qquad \text{for } (t_1, t_2) \in \mathbb{R}^2$$

is easily seen to be stationary and simple max-stable. The a.s. continuity of the process follows from an extension of the arguments in [1].

The specific models we consider are analogous to the ones in the one-dimensional situation:

$$(1.5) \qquad\qquad \phi(t_1, t_2) = \frac{\beta^2}{2\pi} \exp\left\{-\frac{\beta^2(t_1^2 + t_2^2)}{2}\right\}$$

(we call this the normal model),

$$(1.6) \qquad\qquad \phi(t_1, t_2) = \frac{\beta^2}{4} \exp\{-\beta(|t_1| + |t_2|)\}$$

(we call this the exponential model) and

$$(1.7) \qquad \phi(t_1, t_2) = \frac{\beta^2}{2\pi} \left\{1 + \frac{\beta^2(t_1^2 + t_2^2)}{2(\alpha-1)}\right\}^{-\alpha}, \qquad \alpha > 1$$

(we call this the $t$-model), for $\beta > 0$. Finally we consider the general normal model

$$(1.8) \quad \phi(t_1, t_2) = \frac{\beta_1 \beta_2}{2\pi\sqrt{1-\rho^2}} \exp\left\{-\frac{1}{2(1-\rho^2)}[\beta_1^2 t_1^2 - 2\rho\beta_1\beta_2 t_1 t_2 + \beta_2^2 t_2^2]\right\},$$

where $\rho$ is the correlation coefficient $(-1 < \rho < 1)$, for $\beta_1, \beta_2 > 0$.

The paper is organized as follows. The two-dimensional distributions of the process are derived for the mentioned models in Section 2. Those are sufficient for the estimation theory developed in Section 3. The two-dimensional marginal distribution for the normal model was derived earlier by Smith in an unpublished paper [14]. Higher-dimensional marginal distributions do not seem easy to calculate explicitly.

In Section 3 we are mainly concerned with estimating the dependence parameter $\beta$ on the basis of observations at finitely many locations from a stationary stochastic process in the domain of attraction of the max-stable process. In all the models there is a simple relation between $\beta$ and a well-known dependence coefficient for two-dimensional extremes,

$$\lambda := \lim_{t \downarrow 0} t^{-1} P\{1 - F_1(X_1) \leq t \text{ and } 1 - F_2(X_2) \leq t\},$$

where $(X_1, X_2)$ has a distribution function $F$ which is in the domain of attraction of some extreme-value distribution ($F_1$ and $F_2$ are the marginal distribution functions).

The coefficient $\lambda \in [0, 1]$ is related to the general framework of multidimensional extremes in the following way. If

$$\lim_{n \to \infty} F^n(a_n x + b_n, c_n y + d_n) = G(x, y),$$

where the two marginals of $G$ are of the form $\exp\{-(1 + \gamma_i x)_+^{-1/\gamma_i}\}, i = 1, 2$, then

$$\lim_{t \downarrow 0} t^{-1} P\{1 - F_1(X_1) \leq tx \text{ or } 1 - F_2(X_2) \leq ty\}$$

$$= -\log G\left(\frac{x^{-\gamma_1} - 1}{\gamma_1}, \frac{y^{-\gamma_2} - 1}{\gamma_2}\right)$$

$$=: L(x, y)$$

and

$$\lim_{t \downarrow 0} t^{-1} P\{1 - F_1(X_1) \leq tx \text{ and } 1 - F_2(X_2) \leq ty\}$$

$$= x + y - L(x, y) =: R(x, y).$$

Then $\lambda \equiv R(1, 1)$.

D. Mason and X. Huang (see [9]) proved consistency and asymptotic normality for a natural estimator $\tilde{R}(x, y)$ of $R(x, y)$. We use this for estimating $\lambda = R(1, 1)$. This result leads to a consistent asymptotically normal estimator of $\beta$ based on observations taken at just two sites, $t_i$ and $t_j$. In general observations are available at sites $t_1, t_2, \ldots, t_d$ (i.e., finitely many). One of our estimators of $\beta$ is the average of $\beta$-estimators based on the various pairs of sites. Consistency and asymptotic normality follow.

The theory we develop can be used to solve some common problems in spatial extremes. One is an extrapolation problem: it consists of estimating the extremal behavior of a process $X(t)$ at a site $t_0$ where no observations are available based on repeated observations of the process at $d$ different space points $t_1, t_2, \ldots, t_d$. Another problem refers to the extreme values of the (unobserved) aggregate process $\int_S X(t) \, dt$ over a space region $S$ assuming that the process has been discretely observed at a number of space points in $S$. Similarly one can look at the tail behavior of $\sup_{t \in S} X(t)$.

In order to attack those problems, the estimation of $\beta$ has to be complemented by estimation of local parameters: the extreme-value index, the scale and the location. The latter objects are not needed for the estimation of $\beta$ but only for the application.

The proposed models are quite simple examples of the representation in [5], which is valid for all stationary max-stable processes. It seems that the full model is not easily applicable: how does one estimate the initial nonnegative $L_1$ function and a general group of transformations? Moreover the representation is not unique. Instead we have tried to look at simpler models which can be analyzed mathematically. We hope that providing several simple models with quite different features will widen the scope for applications. Later we want to consider somewhat more general parametric groups of transformations.

The validity of the model in applications can be checked with the following steps: estimate $\beta$ as outlined in Section 3, estimate $R(1, 1)$ for each pair of sites using the relation in Corollary 3.2 and check if this estimate of $R(1, 1)$ is similar to the direct estimate using Proposition 3.2. We have not done this yet.

**2. Marginal distributions for the two models, and in $\mathbb{R}$ and $\mathbb{R}^2$.** We find the two-dimensional marginal distributions for all the models introduced in Section 1. Note that in this section, different from other sections, for simplicity and without loss of generality we consider the standardized process $Z$, not $\tilde{Z}$. First we consider the processes on the line, next the ones on $\mathbb{R}^2$.

PROPOSITION 2.1. *For $t \in \mathbb{R}$ and the exponential model, for $w_1, w_2 > 0$,*

$$-\log P\{Z(0) \leq w_1, Z(t) \leq w_2\}$$

$$(2.1) \qquad = \frac{\beta}{2} \int_{-\infty}^{\infty} \max\left\{ \frac{e^{-\beta|s|}}{w_1}, \frac{e^{-\beta|s-t|}}{w_2} \right\} ds$$

$$(2.2) \qquad = \begin{cases} \dfrac{1}{w_2}, & \text{for } 0 < w_2 < e^{-\beta|t|} w_1, \\[2mm] \dfrac{1}{w_1} + \dfrac{1}{w_2} - \dfrac{e^{-\beta|t|/2}}{\sqrt{w_1 w_2}}, & \text{for } e^{-\beta|t|} w_1 \le w_2 < e^{\beta|t|} w_1, \\[2mm] \dfrac{1}{w_1}, & \text{for } w_2 \ge e^{\beta|t|} w_1, \end{cases}$$

$$(2.3) \qquad = \int_{\arctan(e^{-\beta|t|})}^{\arctan(e^{\beta|t|})} \max\left\{ \frac{\sin\theta}{w_1}, \frac{\cos\theta}{w_2} \right\} s(\theta)\, d\theta$$

$$+ \frac{1}{2} \max\left\{ \frac{e^{-\beta|t|}}{w_1}, \frac{1}{w_2} \right\} + \frac{1}{2} \max\left\{ \frac{1}{w_1}, \frac{e^{-\beta|t|}}{w_2} \right\},$$

*where the spectral density $s$ is given by*

$$(2.4) \qquad (\theta) = \frac{e^{-\beta|t|/2}}{4} (\sin\theta \cos\theta)^{-3/2};$$

$$= \frac{1}{w_1} + \frac{1}{w_2} + \frac{1}{w_2} \chi\left( \frac{w_2}{w_1} \right),$$

*and where the dependence function $\chi$ is given by*

$$\chi(s) = \begin{cases} -s, & 0 < s \le e^{-\beta|t|}, \\ -e^{-\beta|t|/2}\sqrt{s}, & e^{-\beta|t|} < s \le e^{\beta|t|}, \\ -1, & s > e^{\beta|t|}, \end{cases} \qquad s > 0.$$

REMARK 2.1.   Formula (2.1) follows directly from [5] and (2.3) reveals the spectral measure, which has a density on the interval $(\arctan(e^{-\beta|t|})$, $\arctan(e^{\beta|t|}))$ and atoms of size $\sqrt{1 + e^{-2\beta|t|}}/2$ at each of the two boundary points of that interval. The last characterization (2.4) is in the spirit of Sibuya [13] and Pickands [11].

PROOF OF PROPOSITION 2.1.   The integrand of (2.1) is $\frac{\beta}{2} \frac{e^{-\beta|s-t|}}{w_2}$ if $\frac{e^{-\beta|s-t|}}{w_2} > \frac{e^{-\beta|s|}}{w_1}$, that is, if

$$(2.5) \qquad |s - t| - |s| < \frac{1}{\beta} \log\left( \frac{w_1}{w_2} \right).$$

Since the joint distributions of $(Z(0), Z(t))$ and $(Z(0), Z(-t))$ are the same we proceed as if $t$ is positive. The left-hand side is $t$ for $s < 0$, $t - 2s$ for $0 < s < t$ and $-t$ for $s > t$. Hence if $\frac{1}{\beta} \log\left( \frac{w_1}{w_2} \right) > t$, then inequality (2.5) holds for all $s$ and we get the first line of (2.2). Similarly, if $\frac{1}{\beta} \log\left( \frac{w_1}{w_2} \right) < -t$, we get the last line of (2.2). Next suppose $-t < \frac{1}{\beta} \log\left(\frac{w_1}{w_2}\right) < t$. In this case (2.5) becomes

$$t - 2s < \frac{1}{\beta} \log\left(\frac{w_1}{w_2}\right),$$

that is,

$$s > \frac{t}{2} - \frac{1}{2\beta} \log\left(\frac{w_1}{w_2}\right).$$

Hence the integrand over this interval becomes

$$\frac{\beta}{2} \int_{t/2 - 1/(2\beta) \log(w_1/w_2)}^{\infty} \frac{e^{-\beta|s-t|}}{w_2} \, ds = \frac{1}{w_2} - \frac{\beta}{2} \int_{-\infty}^{t/2 - 1/(2\beta) \log(w_1/w_2)} \frac{e^{\beta(s-t)}}{w_2} \, ds$$

$$= \frac{1}{w_2} - \frac{1}{2w_2} e^{-\beta t} e^{\beta t/2} e^{-(1/2) \log(w_1/w_2)}$$

$$= \frac{1}{w_2} - \frac{1}{2} \frac{e^{-\beta|t|/2}}{\sqrt{w_1 w_2}}.$$

This in combination with the integral stemming from the case $\frac{e^{-\beta|s-t|}}{w_2} < \frac{e^{-\beta|s|}}{w_1}$ gives the second line of (2.2).

To check the equivalence of (2.2) and (2.3) for $e^{-\beta|t|} < \frac{w_1}{w_2} < e^{\beta|t|}$ it suffices to see that the density of (2.2), after transformation to the polar coordinates $r = \sqrt{w_1^2 + w_2^2}$ and $\theta = \arctan\frac{w_2}{w_1}$, is $r^{-3} s(\theta)$ (cf. the construction of the spectral measure in [6]). For $(w_1, w_2)$ outside this range just evaluate the integral in (2.3). $\square$

REMARK 2.2. The parameter $\beta$ controls the dependence: if $\beta \to \infty$, the spectral density $s(\theta)$ goes to zero and the spectral measure is concentrated at the points $\{0, \frac{\pi}{2}\}$. This means that $X(0)$ and $X(t)$ are independent. If $\beta \downarrow 0$, the spectral measure concentrates on $\{\frac{\pi}{4}\}$. This means that $X(0) = X(t)$ a.s.

PROPOSITION 2.2. *For $t \in \mathbb{R}$ and the normal model, for $w_1, w_2 > 0$,*

$$-\log P\{Z(0) \le w_1, Z(t) \le w_2\}$$

$$(2.6) \qquad = \frac{1}{w_1} \Phi\left(\frac{\beta|t|}{2} + \frac{1}{\beta|t|} \log\frac{w_2}{w_1}\right) + \frac{1}{w_2} \Phi\left(\frac{\beta|t|}{2} + \frac{1}{\beta|t|} \log\frac{w_1}{w_2}\right)$$

$$(2.7) \qquad = \int_0^{\pi/2} \max\left\{\frac{\sin\theta}{w_1}, \frac{\cos\theta}{w_2}\right\} s(\theta) \, d\theta,$$

*where the spectral density $s$ is given by*

$$s(\theta) := \frac{1}{\beta|t| \sin\theta \cos\theta} \left\{\frac{1}{\cos\theta}\left(\frac{1}{2} - \frac{1}{t^2\beta^2} \ln(\tan\theta)\right) \varphi\left(\frac{|t|\beta}{2} + \frac{1}{|t|\beta} \ln(\tan\theta)\right)\right.$$

$$+ \frac{1}{\sin\theta}\left(\frac{1}{2} + \frac{1}{t^2\beta^2}\ln(\tan\theta)\right)\varphi\left(\frac{|t|\beta}{2} - \frac{1}{|t|\beta}\ln(\tan\theta)\right)\Bigg\}$$

*with $\varphi(u) = \Phi'(u)$;*

$$(2.8) \qquad = \frac{1}{w_1} + \frac{1}{w_2} + \frac{1}{w_2}\chi\left(\frac{w_2}{w_1}\right),$$

*where the dependence function $\chi$ is given by*

$$\chi(s) = -s\Phi\left(-\frac{|t|\beta}{2} - \frac{1}{|t|\beta}\ln s\right) - \Phi\left(-\frac{|t|\beta}{2} + \frac{1}{|t|\beta}\ln s\right), \qquad s > 0.$$

REMARK 2.3. Again the parameter $\beta$ controls the dependence: if $\beta \to \infty$, the variables $Z(0)$ and $Z(t)$ become independent; if $\beta \downarrow 0$, we get $Z(0) = Z(t)$ a.s.

REMARK 2.4. This distribution function has been obtained in a number of ways in the literature. Eddy [7] found it when studying convex hulls of samples. Hüsler and Reiss [10] obtained the distribution as the limit distribution of the componentwise maxima in a triangular array where the distribution of the $n$th array is the two-dimensional normal distribution with correlation coefficient $\rho(n)$ such that $\lim_{n\to\infty}(1-\rho(n))\log n$ exists. A related reference is [2]. In [14] Smith developed the distribution in the same way as it is done here. The distribution is mentioned in [3] and [12]. Falk, Hüsler and Reiss [8] obtained the distribution as the pointwise maximum of independent Brownian motions shifted by an amount corresponding to points of a Poisson point process. Another way of obtaining the distribution is

$$-\log P\{Z(0) \le w_1, Z(t) \le w_2\}$$
$$= E\max\left\{\frac{1}{w_1}, \frac{1}{w_2}\exp\left(N\beta t - \beta^2 t^2/2\right)\right\}, \qquad t \in \mathbb{R},$$

where $N$ is a standard normal random variable.

PROOF OF PROPOSITION 2.2. Clearly the distribution depends only on $|t|$. So we consider $t > 0$ only:

$$-\log P\{Z(0) \le w_1, Z(t) \le w_2\} = \frac{\beta}{\sqrt{2\pi}}\int_{-\infty}^{\infty}\max\left\{\frac{e^{-\beta^2 u^2/2}}{w_1}, \frac{e^{-\beta^2(u-t)^2/2}}{w_2}\right\}du.$$

Now $\frac{e^{-\beta^2 u^2/2}}{w_1} \ge \frac{e^{-\beta^2(u-t)^2/2}}{w_2}$ if and only if $\beta u \le \frac{\beta t}{2} + \frac{1}{\beta t}\ln\frac{w_2}{w_1}$, hence

$$-\log P\{Z(0) \le w_1, Z(t) \le w_2\}$$

$$= \frac{1}{w_1} \frac{1}{\sqrt{2\pi}} \int_{-\infty}^{\beta t/2 + (1/\beta t) \ln w_2/w_1} e^{-u^2/2} \, du$$

$$+ \frac{1}{w_2} \frac{1}{\sqrt{2\pi}} \int_{\beta t/2 + (1/\beta t) \ln w_2/w_1}^{\infty} e^{-(u-\beta t)^2/2} \, du$$

$$= \frac{1}{w_1} \Phi\left(\frac{\beta t}{2} + \frac{1}{\beta t} \ln \frac{w_2}{w_1}\right) + \frac{1}{w_2} \left\{ 1 - \Phi\left(-\frac{\beta t}{2} + \frac{1}{\beta t} \ln \frac{w_2}{w_1}\right) \right\}.$$

Hence the first part of the result. In order to obtain the second result note that

$$r^{-3} s(\theta) = \left( -\frac{\partial^2}{\partial w_1 \, \partial w_2} [-\log P\{Z(0) \leq w_1, Z(t) \leq w_2\}] \right)_{w_1 = r\cos\theta, w_2 = r\sin\theta}$$

(cf. the construction of the spectral measure in [6]). □

PROPOSITION 2.3. *For $t \in \mathbb{R}$ and the $t$-model, for $w_1, w_2 > 0$,*

$$-\log P\{Z(0) \leq w_1, Z(t) \leq w_2\}$$

$$(2.9) \qquad = \begin{cases} \dfrac{1}{w_2}, & 0 < w_2 < b_{2,\nu}^{-(\nu+1)/2} w_1, \\[2mm] \dfrac{1}{w_1} p_{1,\nu}(\beta, t, x) + \dfrac{1}{w_2}(1 - p_{2,\nu}(\beta, t, x)), \\[2mm] & b_{2,\nu}^{-(\nu+1)/2} w_1 \leq w_2 < w_1, \\[2mm] \dfrac{2}{w} P\left\{ T_{\nu,1} \leq \dfrac{\beta|t|}{2} \right\}, & w_1 = w_2 =: w, \\[2mm] \dfrac{1}{w_1}(1 - p_{1,\nu}(\beta, t, x)) + \dfrac{1}{w_2} p_{2,\nu}(\beta, t, x), \\[2mm] & w_1 < w_2 < b_{1,\nu}^{-(\nu+1)/2} w_1, \\[2mm] \dfrac{1}{w_1}, & w_2 \geq b_{1,\nu}^{-(\nu+1)/2} w_1, \end{cases}$$

*where*

$$b_{1,\nu} = 1 + \frac{\beta^2 t^2}{2\nu} - \frac{\beta|t|}{\sqrt{\nu}} \sqrt{1 + \frac{\beta^2 t^2}{4\nu}},$$

$$b_{2,\nu} = 1 + \frac{\beta^2 t^2}{2\nu} + \frac{\beta|t|}{\sqrt{\nu}} \sqrt{1 + \frac{\beta^2 t^2}{4\nu}},$$

$$p_{1,\nu}(\beta, t, x) = P\left\{ \left| T_{\nu,1} - \frac{\beta t}{1-x} \right| \leq \beta \sqrt{\frac{t^2 x}{(1-x)^2} - \frac{\nu}{\beta^2}} \right\},$$

$$p_{2,\nu}(\beta, t, x) = P\left\{ \left| T_{\nu,1} - \frac{\beta t x}{1-x} \right| \leq \beta \sqrt{\frac{t^2 x}{(1-x)^2} - \frac{\nu}{\beta^2}} \right\},$$

$T_{\nu,1}$ *is a random variable with a Student t-distribution with $\nu$ degrees of freedom and scale parameter 1 and* $x = (\frac{w_1}{w_2})^{2/(\nu+1)}$;

$$(2.10) \qquad\qquad = \frac{1}{w_1} + \frac{1}{w_2} + \frac{1}{w_2}\chi\left(\frac{w_2}{w_1}\right),$$

*where the dependence function $\chi$ is given by*

$$\chi(s) = \begin{cases} -s, & 0 < s \leq b_{2,\nu}^{-(\nu+1)/2}, \\ sp_{1,\nu}(\beta, t, s^{-2/(\nu+1)}) - s - p_{2,\nu}(\beta, t, s^{-2/(\nu+1)}), & \\ & b_{2,\nu}^{-(\nu+1)/2} < s \leq 1, \\ -sp_{1,\nu}(\beta, t, s^{-2/(\nu+1)}) - 1 + p_{2,\nu}(\beta, t, s^{-2/(\nu+1)}), & \\ & 1 < s \leq b_{1,\nu}^{-(\nu+1)/2}, \\ -1, & s > b_{1,\nu}^{-(\nu+1)/2}, \end{cases}$$

$s > 0$.

REMARK 2.5.  The spectral measure corresponding to (2.9) is concentrated on $[b_{2,\nu}^{-(\nu+1)/2}, b_{1,\nu}^{-(\nu+1)/2}]$, having a density on $(b_{2,\nu}^{-(\nu+1)/2}, b_{1,\nu}^{-(\nu+1)/2})$ and atoms at each of the boundary points of the interval.

PROOF OF PROPOSITION 2.3.

$-\log P\{Z(0) \leq w_1, Z(t) \leq w_2\}$

$$= \frac{\beta\Gamma((\nu+1)/2)}{\sqrt{\pi\nu}\Gamma(\nu/2)}$$

$$\times \int_{-\infty}^{\infty} \max\left\{ \frac{(1+\beta^2 u^2/\nu)^{-(\nu+1)/2}}{w_1}, \frac{(1+\beta^2(u-t)^2/\nu)^{-(\nu+1)/2}}{w_2} \right\} du.$$

Now $\frac{(1+\beta^2 u^2/\nu)^{-(\nu+1)/2}}{w_1} \geq \frac{(1+\beta^2(u-t)^2/\nu)^{-(\nu+1)/2}}{w_2}$ if and only if $(1-x)u^2 - 2tu \geq -\frac{(1-x)\nu}{\beta^2} - t^2$. This is equivalent to

$$(2.11) \qquad\qquad \left(u - \frac{t}{1-x}\right)^2 \geq \frac{xt^2}{(1-x)^2} - \frac{\nu}{\beta^2}$$

if $0 < x < 1$ and to the reversed inequality if $x > 1$. Hence if $0 < x \leq b_{1,\nu}$, then (2.11) holds for all $u$ and we get the last line of (2.9). Similarly if $x > b_{2,\nu}$, we get the first line of (2.9). In the case $b_{1,\nu} < x < 1$ condition (2.11) becomes

$$\left| u - \frac{t}{1-x} \right| \geq \sqrt{\frac{xt^2}{(1-x)^2} - \frac{\nu}{\beta^2}},$$

that is,

$$u \geq \frac{t}{1-x} + \sqrt{\frac{xt^2}{(1-x)^2} - \frac{\nu}{\beta^2}} \quad \text{or} \quad u \leq \frac{t}{1-x} - \sqrt{\frac{xt^2}{(1-x)^2} - \frac{\nu}{\beta^2}}.$$

Hence the integral over this interval becomes

$$\frac{\beta\Gamma((\nu+1)/2)}{\sqrt{\pi\nu}\Gamma(\nu/2)} \int_{-\infty}^{t/(1-x)-\sqrt{xt^2/(1-x)^2-\nu/\beta^2}} \frac{(1+\beta^2 u^2/\nu)^{-(\nu+1)/2}}{w_1} du$$

$$+ \frac{\beta\Gamma((\nu+1)/2)}{\sqrt{\pi\nu}\Gamma(\nu/2)} \int_{t/(1-x)+\sqrt{xt^2/(1-x)^2-\nu/\beta^2}}^{+\infty} \frac{(1+\beta^2 u^2/\nu)^{-(\nu+1)/2}}{w_1} du$$

$$= \frac{1}{w_1} \left\{ F_{\nu,1}\left( \frac{\beta t}{1-x} - \beta\sqrt{\frac{xt^2}{(1-x)^2} - \frac{\nu}{\beta^2}} \right) \right.$$

$$\left. + 1 - F_{\nu,1}\left( \frac{\beta t}{1-x} + \beta\sqrt{\frac{xt^2}{(1-x)^2} - \frac{\nu}{\beta^2}} \right) \right\}.$$

This, in combination with the integral stemming from the case $\frac{(1+\beta^2 u^2/\nu)^{-(\nu+1)/2}}{w_1} < \frac{(1+\beta^2(u-t)^2/\nu)^{-(\nu+1)/2}}{w_2}$, gives the fourth line of (2.9). Similarly in the case $1 < x \leq b_{2,\nu}$ we get the second line of (2.9). Note that in the case $x = 1$, $\frac{(1+\beta^2 u^2/\nu)^{-(\nu+1)/2}}{w_1} \geq \frac{(1+\beta^2(u-t)^2/\nu)^{-(\nu+1)/2}}{w_2}$ if and only if $u \leq t/2$ for $t > 0$ and $u \geq t/2$ for $t < 0$. Hence we get for (2.9) $\frac{2}{w_1} P\{T_{\nu,1} \leq \beta|t|/2\}$. $\quad \square$

Next we move to the two-dimensional models.

PROPOSITION 2.4. *For $t = (t_1, t_2) \in \mathbb{R}^2$ and the exponential model, for $w_1, w_2 > 0$,*

$$-\log P\{Z(0,0) \leq w_1, Z(t_1,t_2) \leq w_2\}$$

$$(2.12) \quad = \begin{cases} \dfrac{1}{w_1}, & (w_1, w_2) \in A_1, \\[2mm] \dfrac{1}{w_1} + \dfrac{1}{w_2} - \dfrac{1}{\sqrt{w_1 w_2}} e^{-\beta(|t_1|+|t_2|)/2} \left[ 1 + \beta\dfrac{|t_1|+|t_2|}{4} + \dfrac{1}{4}\ln\left(\dfrac{w_1}{w_2}\right) \right], \\ & (w_1, w_2) \in A_2, \\[2mm] \dfrac{1}{w_1} + \dfrac{1}{w_2} - \dfrac{1}{\sqrt{w_1 w_2}} e^{-\beta(|t_1|+|t_2|)/2} \left[ 1 + \beta\dfrac{\min(|t_1|,|t_2|)}{2} \right], \\ & (w_1, w_2) \in A_3, \\[2mm] \dfrac{1}{w_1} + \dfrac{1}{w_2} - \dfrac{1}{\sqrt{w_1 w_2}} e^{-\beta(|t_1|+|t_2|)/2} \left[ 1 + \beta\dfrac{|t_1|+|t_2|}{4} + \dfrac{1}{4}\ln\left(\dfrac{w_2}{w_1}\right) \right], \\ & (w_1, w_2) \in A_4, \\[2mm] \dfrac{1}{w_2}, & (w_1, w_2) \in A_5, \end{cases}$$

*where*

$$A_1 = \left\{ (w_1, w_2) : \frac{1}{\beta} \ln\left(\frac{w_1}{w_2}\right) < -(|t_1| + |t_2|) \right\},$$

$$A_2 = \left\{ (w_1, w_2) : -(|t_1| + |t_2|) \leq \frac{1}{\beta} \ln\left(\frac{w_1}{w_2}\right) < |t_1| \wedge |t_2| - |t_1| \vee |t_2| \right\},$$

$$A_3 = \left\{ (w_1, w_2) : |t_1| \wedge |t_2| - |t_1| \vee |t_2| \leq \frac{1}{\beta} \ln\left(\frac{w_1}{w_2}\right) < |t_1| \vee |t_2| - |t_1| \wedge |t_2| \right\},$$

$$A_4 = \left\{ (w_1, w_2) : |t_1| \vee |t_2| - |t_1| \wedge |t_2| \leq \frac{1}{\beta} \ln\left(\frac{w_1}{w_2}\right) < |t_1| + |t_2| \right\},$$

$$A_5 = \left\{ (w_1, w_2) : \frac{1}{\beta} \ln\left(\frac{w_1}{w_2}\right) \geq |t_1| + |t_2| \right\}.$$

Proof. We work out the integral as in the one-dimensional case on the areas of $\mathbf{u}$ defined by $|u_1 - t_1| + |u_2 - t_2| - |u_1| - |u_2| \gtrless \frac{1}{\beta} \ln(\frac{w_1}{w_2})$. It is convenient, similarly to the one-dimensional case, to consider separately the nine areas defined by the position of $u_1$ with respect to 0 and $t_1$ and by the position of $u_2$ with respect to 0 and $t_2$. So, for example, if $0 < t_1 < t_2$, we can write

$$|u_1 - t_1| + |u_2 - t_2| - |u_1| - |u_2|$$

$$= \begin{cases} -t_1 - t_2, & \text{if } u_1 > t_1 \text{ and } u_2 > t_2, \\ -2u_2 + t_2 - t_1, & \text{if } u_1 > t_1 \text{ and } 0 < u_2 < t_2, \\ t_2 - t_1, & \text{if } u_1 > t_1 \text{ and } u_2 < 0, \\ -2u_1 + t_1 - t_2, & \text{if } 0 < u_1 < t_1 \text{ and } u_2 > t_2, \\ -2u_1 - 2u_2 + t_1 + t_2, & \text{if } 0 < u_1 < t_1 \text{ and } 0 < u_2 < t_2, \\ -2u_1 + t_1 + t_2, & \text{if } 0 < u_1 < t_1 \text{ and } u_2 < 0, \\ t_1 - t_2, & \text{if } u_1 < 0 \text{ and } u_2 > t_2, \\ -2u_2 + t_2 + t_1, & \text{if } u_1 < 0 \text{ and } 0 < u_2 < t_2, \\ t_1 + t_2, & \text{if } u_1 < 0 \text{ and } u_2 < 0. \end{cases}$$

The calculations are complicated but not difficult. □

PROPOSITION 2.5. *For* $t = (t_1, t_2) \in \mathbb{R}^2$ *and the normal model, we have for* $w_1, w_2 > 0$,

$$(2.13) \quad \begin{aligned} &-\log P\{Z(0,0) \leq w_1, Z(t_1, t_2) \leq w_2\} \\ &= \frac{1}{w_1} \Phi\left(\frac{|\mathbf{t}|\beta}{2} + \frac{1}{|\mathbf{t}|\beta} \log \frac{w_2}{w_1}\right) + \frac{1}{w_2} \Phi\left(\frac{|\mathbf{t}|\beta}{2} + \frac{1}{|\mathbf{t}|\beta} \log \frac{w_1}{w_2}\right), \end{aligned}$$

*that is, it is the same as in the one-dimensional case with* $|t|$ *replaced by* $|\mathbf{t}|$.

PROOF. As in the one-dimensional case. Now $\frac{e^{-\beta^2|\mathbf{u}|^2/2}}{w_1} \geq \frac{e^{-\beta^2|\mathbf{u}-\mathbf{t}|^2/2}}{w_2}$ if and only if $\beta \frac{t_1 u_1 + t_2 u_2}{|\mathbf{t}|} \leq \frac{|\mathbf{t}|\beta}{2} + \frac{1}{|\mathbf{t}|\beta} \log \frac{w_2}{w_1}$. Next note that if the vector $(U_1, U_2)$ has a standard two-dimensional normal distribution, $\frac{t_1 U_1 + t_2 U_2}{|\mathbf{t}|}$ is also standard normal. The rest of the proof is as in the one-dimensional case. $\square$

PROPOSITION 2.6. *For $t = (t_1, t_2) \in \mathbb{R}^2$ and the $t$-model, for $w_1, w_2 > 0$,*

$$-\log P\{Z(0,0) \leq w_1, Z(t_1, t_2) \leq w_2\}$$

$$(2.14) \quad = \begin{cases} \dfrac{1}{w_2}, & 0 < w_2 < b_{2,\alpha}^{-\alpha} w_1, \\[2mm] \dfrac{1}{w_1} P\{(T_1, T_2) \in A_{1,\alpha}\} + \dfrac{1}{w_2} P\{(T_1, T_2) \in A_{2,\alpha}^c\}, \\[1mm] & b_{2,\alpha}^{-\alpha} w_1 \leq w_2 < w_1, \\[2mm] \dfrac{2}{w} P\{T_1 \leq |\mathbf{t}|/2\}, & w_1 = w_2 =: w, \\[2mm] \dfrac{1}{w_1} P\{(T_1, T_2) \in A_{1,\alpha}^c\} + \dfrac{1}{w_2} P\{(T_1, T_2) \in A_{2,\alpha}\}, \\[1mm] & w_1 < w_2 < b_{1,\alpha}^{-\alpha} w_1, \\[2mm] \dfrac{1}{w_1}, & w_2 \geq b_{1,\alpha}^{-\alpha} w_1, \end{cases}$$

*where $(T_1, T_2)$ is a random vector with bivariate $t$-density* (1.7)*,*

$$b_{1,\alpha} = 1 + \frac{\beta^2|\mathbf{t}|^2}{4(\alpha-1)} - \frac{\beta|\mathbf{t}|}{\sqrt{2(\alpha-1)}}\sqrt{1 + \frac{\beta^2|\mathbf{t}|^2}{8(\alpha-1)}},$$

$$b_{2,\alpha} = 1 + \frac{\beta^2|\mathbf{t}|^2}{4(\alpha-1)} + \frac{\beta|\mathbf{t}|}{\sqrt{2(\alpha-1)}}\sqrt{1 + \frac{\beta^2|\mathbf{t}|^2}{8(\alpha-1)}},$$

$$A_{1,\alpha} = \left\{ (u_1, u_2) \in \mathbb{R}^2 : \left(u_1 - \frac{t_1}{1-x}\right)^2 + \left(u_2 - \frac{t_2}{1-x}\right)^2 \right.$$
$$\left. \leq \frac{x|\mathbf{t}|^2}{(1-x)^2} - \frac{2(\alpha-1)}{\beta^2} \right\},$$

$$A_{2,\alpha} = \left\{ (u_1, u_2) \in \mathbb{R}^2 : \left(u_1 - \frac{t_1 x}{1-x}\right)^2 + \left(u_2 - \frac{t_2 x}{1-x}\right)^2 \right.$$
$$\left. \leq \frac{x|\mathbf{t}|^2}{(1-x)^2} - \frac{2(\alpha-1)}{\beta^2} \right\}$$

*and*

$$x = \left(\frac{w_1}{w_2}\right)^{1/\alpha}.$$

PROOF.    Analogous to the one-dimensional case.    □

PROPOSITION 2.7.    *For* $t = (t_1, t_2) \in \mathbb{R}^2$ *and the general normal model, for* $w_1, w_2 > 0$,

$$-\log P\{Z(0,0) \leq w_1, Z(t_1, t_2) \leq w_2\}$$

$$(2.15) \qquad = \frac{1}{w_1} \Phi\left( \frac{\sqrt{\mathbf{t}^T \Sigma^{-1} \mathbf{t}}}{2} + \frac{1}{\sqrt{\mathbf{t}^T \Sigma^{-1} \mathbf{t}}} \log \frac{w_2}{w_1} \right)$$

$$+ \frac{1}{w_2} \Phi\left( \frac{\sqrt{\mathbf{t}^T \Sigma^{-1} \mathbf{t}}}{2} + \frac{1}{\sqrt{\mathbf{t}^T \Sigma^{-1} \mathbf{t}}} \log \frac{w_1}{w_2} \right),$$

*where*

$$\Sigma^{-1} = \frac{1}{1 - \rho^2} \begin{bmatrix} \beta_1^2 & -\rho \beta_1 \beta_2 \\ -\rho \beta_1 \beta_2 & \beta_2^2 \end{bmatrix}.$$

PROOF.    As in the one-dimensional case. Now $\frac{e^{-\mathbf{u}^T \Sigma^{-1} \mathbf{u}/2}}{w_1} \geq \frac{e^{-(\mathbf{u} - \mathbf{t})^T \Sigma^{-1} (\mathbf{u} - \mathbf{t})/2}}{w_2}$ if and only if $\mathbf{u}^T \Sigma^{-1} \mathbf{t} \leq \frac{\mathbf{t}^T \Sigma^{-1} \mathbf{t}}{2} + \log \frac{w_2}{w_1}$. Next note that if the vector $\mathbf{U} = (U_1, U_2)$ has a bivariate normal distribution with mean value zero and covariance matrix $\Sigma$, $\mathbf{U}^T \Sigma^{-1} \mathbf{t}$ has a normal distribution with mean value zero and variance $\mathbf{t}^T \Sigma^{-1} \mathbf{t}$. The rest of the proof is as in the one-dimensional case.    □

**3. Estimating the dependence parameter** $\boldsymbol{\beta}$**.**    We consider a sequence of independent, identically distributed stochastic processes with continuous paths

$$\{X_i(t)\}_{t \in \mathbb{R}}, \qquad i = 1, 2, \ldots.$$

We assume that the processes are in the max-domain of attraction [as processes in $C(\mathbb{R})$] of a max-stable stationary process $\{\tilde{Z}(t)\}_{t \in \mathbb{R}}$ such that the related process $Z$ (see the Introduction) has exponential spectral function (2.2) or (2.12), or normal spectral function (2.6) or (2.13), or $t$ spectral function (2.9) or (2.14) as discussed in Section 2. For definition of convergence and convergence criteria see [4].

Ideally it would be nice if we could assume that we have observed the sample paths of $n$ processes $X$ as a basis for estimation of the main parameter $\beta$. However, in reality this is too much to expect. Usually one can observe the $n$ processes only at finitely many points in space, say $t_1, t_2, \ldots, t_d$.

In this setup we propose estimators for $\beta$ that are closely related to an extension of the estimator $\hat{R}(x, y)$ for the dependence function $R(x, y)$ which was introduced by Mason and Huang (see Huang [9]),

$$\hat{R}_{t_1, \ldots, t_d}(x_1, \ldots, x_d) := \frac{1}{k} \sum_{i=1}^{n} I_{\{X_i(t_1) \geq X_{n-[kx_1]+1,n}(t_1), \ldots, X_i(t_d) \geq X_{n-[kx_d]+1,n}(t_d)\}}$$

with $k \ll n$, where $\{X_{i,n}(t_j)\}_{i=1}^n$ are the $n$th-order statistics of $\{X_i(t_j)\}_{i=1}^n$ for $j = 1, 2, \ldots, d$.

Let us now introduce our estimators in $\mathbb{R}$. All three densities have the form $\beta\varphi_0(\beta t)$. Corollary 3.2 below states that

$$R_{t_1,\ldots,t_d}(1,\ldots,1) = 2 \int_{(\beta/2)(\max_{1\le j\le d} t_j - \min_{1\le j\le d} t_j)}^{+\infty} \varphi_0(t)\, dt$$
$$= 2\left\{1 - F\left(\frac{\beta}{2}\left(\max_{1\le j\le d} t_j - \min_{1\le j\le d} t_j\right)\right)\right\}$$

with $F(t) := \int_{-\infty}^t \varphi_0(s)\, ds$. It follows that

$$\beta = \frac{2F^\leftarrow(1 - R_{t_1,t_2,\ldots,t_d}(1,1,\ldots,1)/2)}{\max_{1\le j\le d} t_j - \min_{1\le j\le d} t_j}.$$

Hence we introduce the estimators

$$\hat{\beta}^{(1)} := \frac{2F^\leftarrow(1 - \hat{R}_{t_1,t_2,\ldots,t_d}(1,1,\ldots,1)/2)}{\max_{1\le j\le d} t_j - \min_{1\le j\le d} t_j}$$

and

$$\hat{\beta}_1 := \frac{2}{d(d-1)} \sum_{1\le j<m\le d} \frac{2}{|t_j - t_m|} F^\leftarrow(1 - \hat{R}_{t_j,t_m}(1,1)/2).$$

Considering the two-dimensional space, note that the standard normal and Student densities are spherical symmetric. This means (cf. proof of Proposition 2.5) that it is sufficient to consider the marginal distribution and hence we introduce the estimators

$$\hat{\beta}_2 := \frac{2}{d(d-1)} \sum_{1\le j<m\le d} \frac{2}{|\mathbf{t}_j - \mathbf{t}_m|} F^\leftarrow(1 - \hat{R}_{\mathbf{t}_j,\mathbf{t}_m}(1,1)/2),$$

with $F$ the (common) marginal distribution corresponding to $\varphi_0$ (standard normal or Student). The exponential model in two-dimensional space is more complicated. We consider

$$\hat{\beta}_{e,2} := \frac{2}{d(d-1)} \sum_{1\le j<m\le d} Q_{a_{j,m},b_{j,m}}(\hat{R}_{\mathbf{t}_j,\mathbf{t}_m}(1,1))$$

with $a_{j,m} := |t_1^{(j)} - t_1^{(m)}|$, the absolute difference of the first components of $\mathbf{t}_j$ and $\mathbf{t}_m$, and $b_{j,m} := |t_2^{(j)} - t_2^{(m)}|$ and $Q_{a,b}$ the inverse function of $\frac{1}{2}(1 + \frac{\beta}{2}\min(a,b))\exp\{-\frac{\beta}{2}(a+b)\}$, which is decreasing in $\beta$ for $a, b > 0$.

Note that $\hat{\beta}^{(1)}$ is simpler than $\hat{\beta}$ and summarizes the information of the sample in a somewhat more crude way. We could not find analogues of $\hat{\beta}^{(1)}$ in two-dimensional space since we were unable to calculate explicitly the necessary higher-dimensional distributions.

All the mentioned estimators are consistent and asymptotically normal under the appropriate conditions. We now state the results. The proofs will be given at the end of the section after some lemmas. First we consider consistency.

THEOREM 3.1.   *Suppose that the normalized sequence of maxima [see (1.1)] converges weakly to $\{\tilde{Z}(t)\}$ in $C(-\infty, \infty)$. For sequences $k = k(n) \to \infty$, $k(n)/n \to 0$ as $n \to \infty$ (recall that the sequence $k$ figures in the definition of $\hat{R}$), all the indicated estimators are weakly consistent for $\beta$.*

Also the estimators are asymptotically normal under some extra conditions, that is, $\sqrt{k}(\hat{\beta} - \beta)$ has asymptotically a normal mean zero distribution. In order to describe the asymptotic distribution more accurately, we now state a slight extension of a result of Huang and Mason.

PROPOSITION 3.2 ([9], pages 29 and 43).   *Let $\{(X_i(t_1), \ldots, X_i(t_d))\}_{i=1}^{\infty}$ be i.i.d. random vectors with distribution function $F$. Suppose that the marginal distributions $F_i$ $(i = 1, 2, \ldots, d)$ are continuous and strictly increasing. Define*

$$\tilde{F}_{t_1,\ldots,t_d}(x_1,\ldots,x_d) := 1 - F(F_1^{\leftarrow}(1 - x_1), \ldots, F_d^{\leftarrow}(1 - x_d))$$
$$= P\{1 - F_1(X(t_1)) \leq x_1 \ or \ \cdots \ or \ 1 - F_d(X(t_d)) \leq x_d\},$$

*where the arrow denotes inverse function. Suppose that for all $x_1, x_2, \ldots, x_d \geq 0$, $x_1 + x_2 + \cdots + x_d > 0$ and a positive function $L$,*

$$(3.1) \qquad \lim_{t \downarrow 0} \frac{1}{t} \tilde{F}_{t_1, t_2, \ldots, t_d}(tx_1, tx_2, \ldots, tx_d) = L_{t_1, t_2, \ldots, t_d}(x_1, x_2, \ldots, x_d).$$

*Next we introduce the definition of a function $R$ which is connected with the function $L$ as follows. Let $\nu_{t_1, t_2, \ldots, t_d}$ be the measure that satisfies for $x_1, x_2, \ldots, x_d > 0$*

$$\nu_{t_1,\ldots,t_d}\{(s_1,\ldots,s_d)|s_1 \leq x_1 \ or \ \cdots \ or \ s_d \leq x_d\} = L_{t_1,\ldots,t_d}(x_1,\ldots,x_d).$$

*Such a measure exists by virtue of (3.1). The function $R$ is given by*

$$R_{t_1, t_2, \ldots, t_d}(x_1, x_2, \ldots, x_d) := \nu_{t_1, t_2, \ldots, t_d}([0, x_1] \times [0, x_2] \times \cdots \times [0, x_d]).$$

*Define*

$$\hat{R}_{t_1,\ldots,t_d}(x_1,\ldots,x_d) := \frac{1}{k} \sum_{i=1}^{n} I_{\{X_i(t_1) \geq X_{n-[kx_1]+1,n}(t_1), \ldots, X_i(t_d) \geq X_{n-[kx_d]+1,n}(t_d)\}},$$

*where $X_{n-[kx_i]+1,n}(t_j)$ is the $(n-[kx_i]+1)$st-order statistic of $X_1(t_j), X_2(t_j), \ldots, X_n(t_j), j = 1, 2, \ldots, d$. Then for all $x_1, x_2, \ldots, x_d \geq 0$, $x_1 + x_2 + \cdots + x_d > 0$,*

$$(3.2) \quad \hat{R}_{t_1, t_2, \ldots, t_d}(x_1, x_2, \ldots, x_d) \to R_{t_1, t_2, \ldots, t_d}(x_1, x_2, \ldots, x_d) \qquad in \ probability,$$

$n \to \infty$, $k = k(n) \to \infty$, $k/n \to 0$.

*Suppose now that the function $L$ has continuous first-order partial derivatives $L^{(j)}_{t_1, t_2, \ldots, t_d}$, $j = 1, 2, \ldots, d$. Moreover suppose that for some $\alpha > 0$*

$$(3.3) \qquad \begin{aligned} &\tilde{F}_{t_1, t_2, \ldots, t_d}(tx_1, tx_2, \ldots, tx_d) \\ &\qquad = t\{L_{t_1, t_2, \ldots, t_d}(x_1, x_2, \ldots, x_d) + O(t^\alpha)\}, \qquad t \downarrow 0, \end{aligned}$$

*uniformly on $x_1^2 + x_2^2 + \cdots + x_d^2 = 1$, $x_i \geq 0$, $i = 1, 2, \ldots, d$. Then for a sequence $k = k(n) \to \infty$ with $k = o(n^{2\alpha/(2\alpha+1)})$, $n \to \infty$,*

$$(3.4) \quad \sqrt{k}\{\hat{R}_{t_1, \ldots, t_d}(x_1, \ldots, x_d) - R_{t_1, \ldots, t_d}(x_1, \ldots, x_d)\} \xrightarrow{d} B_{t_1, \ldots, t_d}(x_1, \ldots, x_d)$$

*in $D(\mathbb{R}^d_+)$, where*

$$\begin{aligned} &B_{t_1, t_2, \ldots, t_d}(x_1, x_2, \ldots, x_d) \\ &\quad = W_{t_1, t_2, \ldots, t_d}(x_1, x_2, \ldots, x_d) \\ &\qquad - L^{(1)}_{t_1, t_2, \ldots, t_d}(x_1, x_2, \ldots, x_d)W_{t_1, t_2, \ldots, t_d}(x_1, 0, \ldots, 0) \\ &\qquad - \cdots - L^{(d)}_{t_1, t_2, \ldots, t_d}(x_1, x_2, \ldots, x_d)W_{t_1, t_2, \ldots, t_d}(0, 0, \ldots, x_d) \end{aligned}$$

*and $W$ is a continuous mean zero Gaussian process with covariance structure*

$$EW_{t_1, \ldots, t_d}(x_1, \ldots, x_d)W_{t_1, \ldots, t_d}(y_1, \ldots, y_d) = \nu_{t_1, \ldots, t_d}(A_{x_1, \ldots, x_d} \cap A_{y_1, \ldots, y_d})$$

*with*

$$A_{x_1, x_2, \ldots, x_d} := \{(t_1, t_2, \ldots, t_d) | t_1 < x_1 \text{ or } t_2 < x_2 \text{ or } \cdots \text{ or } t_d < x_d\}.$$

The function $L$ is called the tail dependence function and it is directly related with the extreme-value limit distribution. In fact if $F$ is in the domain of attraction of an extreme-value distribution $G$, condition (3.1) holds with $L(x_1, \ldots, x_d) = -\log G((-\log G_1)^{\leftarrow}(x_1), \ldots, (-\log G_d)^{\leftarrow}(x_d))$ where $G_1, \ldots, G_d$ are the marginals of $G$.

REMARK 3.1.  For the exponential model in one-dimensional space

$$L_{t_1, \ldots, t_d}(x_1, \ldots, x_d) := \frac{\beta}{2}\int_{-\infty}^{\infty} \max(x_1 e^{-\beta|s - t_1|}, \ldots, x_d e^{-\beta|s - t_d|})\, ds$$

and (see Lemma 3.1 below)

$$R_{t_1, \ldots, t_d}(x_1, \ldots, x_d) := \frac{\beta}{2}\int_{-\infty}^{\infty} \min(x_1 e^{-\beta|s - t_1|}, \ldots, x_d e^{-\beta|s - t_d|})\, ds.$$

Similarly for the normal model and the $t$-model.

THEOREM 3.3 (One-dimensional space). *Suppose for the processes* $\{X_n(t)\}$ *the limit relation* (1.1) *holds and moreover the second-order condition* (3.3): *for some* $\alpha > 0$

$$\tilde{F}_{t_1,\dots,t_d}(tx_1,\dots,tx_d) = t\{-\log P\{Z(t_1) \le 1/x_1,\dots,Z(t_d) \le 1/x_d\} + O(t^\alpha)\},$$

$t \downarrow 0$, *uniformly on* $x_1^2 + x_2^2 + \cdots + x_d^2 = 1$, $x_i \ge 0, i = 1,2,\dots,d$.

*If* $k = k(n) \to \infty$, $k(n) = o(n^{2\alpha/(2\alpha+1)})$, *as* $n \to \infty$, *then*

$$\sqrt{k}(\hat{\beta}_1 - \beta) \to \frac{2}{d(d-1)} \sum_{j<m} \frac{2B_{t_j,t_m}(1,1)}{|t_j - t_m|} \frac{1}{\varphi_0(\beta|t_j - t_m|/2)} \tag{3.5}$$

*and*

$$\sqrt{k}(\hat{\beta}^{(1)} - \beta)$$
$$\to \frac{2B_{t_1,t_2,\dots,t_d}(1,1,\dots,1)}{\max_{1\le j\le d} t_j - \min_{1\le j\le d} t_j} \frac{1}{\varphi_0(\beta(\max_{1\le j\le d} t_j - \min_{1\le j\le d} t_j)/2)} \tag{3.6}$$

*in distribution. Here B is as in Proposition* 3.2.

THEOREM 3.4 (Two-dimensional space). *Under the same conditions,*

$$\sqrt{k}(\hat{\beta}_2 - \beta) \to \frac{2}{d(d-1)} \sum_{j<m} \frac{2B_{\mathbf{t}_j,\mathbf{t}_m}(1,1)}{|\mathbf{t}_j - \mathbf{t}_m|} \frac{1}{\varphi_0(\beta|\mathbf{t}_j - \mathbf{t}_m|/2)} \tag{3.7}$$

*and*

$$\sqrt{k}(\hat{\beta}_{e,2} - \beta) \to \frac{2}{d(d-1)} \sum_{j<m} B_{\mathbf{t}_j,\mathbf{t}_m}(1,1) Q'_{a_{jm},b_{jm}}(Q^{\leftarrow}_{a_{jm},b_{jm}}(\beta)) \tag{3.8}$$

*in distribution, where* $a_{jm} := |t_1^{(j)} - t_1^{(m)}|$, $b_{jm} := |t_2^{(j)} - t_2^{(m)}|$.

For the estimation of the general normal model we proceed as follows. Write

$$\Sigma^{-1} = \frac{1}{1-\rho^2} \begin{bmatrix} \beta_1^2 & -\rho\beta_1\beta_2 \\ -\rho\beta_1\beta_2 & \beta_2^2 \end{bmatrix}$$

with $-1 < \rho < 1$ and $\beta_1, \beta_2 > 0$. For two sites $\mathbf{t}_j$ and $\mathbf{t}_m$ $(1 \le j < m \le d)$ we write

$$(\mathbf{t}_j - \mathbf{t}_m)^T \Sigma^{-1} (\mathbf{t}_j - \mathbf{t}_m)$$
$$= \frac{1}{1-\rho^2}\{\beta_1^2(t_j^{(1)} - t_m^{(1)})^2$$
$$- 2\rho\beta_1\beta_2(t_j^{(1)} - t_m^{(1)})(t_j^{(2)} - t_m^{(2)}) + \beta_2^2(t_j^{(2)} - t_m^{(2)})^2\}$$
$$= \mathbf{t}_{j,m}^T \mathbf{a},$$

where

$$\mathbf{t}_j = (t_j^{(1)}, t_j^{(2)})^T \qquad \text{for } j = 1, 2, \ldots, d,$$

$$\mathbf{t}_{j,m} = \begin{bmatrix} (t_j^{(1)} - t_m^{(1)})^2 \\ (t_j^{(1)} - t_m^{(1)})(t_j^{(2)} - t_m^{(2)}) \\ (t_j^{(2)} - t_m^{(2)})^2 \end{bmatrix}$$

and

$$(3.9) \qquad \mathbf{a} := \frac{1}{1-\rho^2} \begin{bmatrix} \beta_1^2 \\ -2\rho\beta_1\beta_2 \\ \beta_2^2 \end{bmatrix}.$$

Now note that

$$2 - R_{\mathbf{t}_j, \mathbf{t}_m}(1,1) = L_{\mathbf{t}_j, \mathbf{t}_m}(1,1) = 2\Phi(\sqrt{|\mathbf{t}_j - \mathbf{t}_m|^T \Sigma^{-1} |\mathbf{t}_j - \mathbf{t}_m|}/2).$$

Hence we define estimators

$$(3.10) \qquad \hat{Q}_{j,m} := (2\Phi^{\leftarrow}(1 - R_{\mathbf{t}_j, \mathbf{t}_m}(1,1)/2))^2.$$

Using the result of Proposition 3.2 and Cramér's delta method we get

$$(3.11) \qquad \sqrt{k}(\hat{Q}_{j,m} - \mathbf{t}_{j,m}^T \mathbf{a}) \xrightarrow{d} 2B_{\mathbf{t}_j, \mathbf{t}_m}(1,1)\sqrt{\mathbf{t}_{j,m}^T \mathbf{a}}(\phi(\sqrt{\mathbf{t}_{j,m}^T \mathbf{a}}/2))^{-1},$$

where $\phi$ is the standard normal density. Now compose the $(d(d-1)/2)$-dimensional vectors

$$\hat{\mathbf{q}} := \begin{bmatrix} \hat{Q}_{1,2} \\ \hat{Q}_{1,3} \\ \vdots \\ \hat{Q}_{d-1,d} \end{bmatrix}$$

and

$$\boldsymbol{\Gamma} := \begin{bmatrix} \mathbf{t}_{1,2}^T \\ \mathbf{t}_{1,3}^T \\ \vdots \\ \mathbf{t}_{d-1,d}^T \end{bmatrix}.$$

Then

$$(3.12) \qquad \sqrt{k}(\hat{\mathbf{q}} - \boldsymbol{\Gamma}\mathbf{a}) \xrightarrow{d} \mathbf{b}$$

with

$$\mathbf{b} := 2 \begin{bmatrix} B_{\mathbf{t}_1,\mathbf{t}_2}(1,1)\sqrt{\mathbf{t}_{1,2}^T\mathbf{a}}(\phi(\sqrt{\mathbf{t}_{1,2}^T\mathbf{a}}/2))^{-1} \\ B_{\mathbf{t}_1,\mathbf{t}_3}(1,1)\sqrt{\mathbf{t}_{1,3}^T\mathbf{a}}(\phi(\sqrt{\mathbf{t}_{1,3}^T\mathbf{a}}/2))^{-1} \\ \vdots \\ B_{\mathbf{t}_{d-1},\mathbf{t}_d}(1,1)\sqrt{\mathbf{t}_{d-1,d}^T\mathbf{a}}(\phi(\sqrt{\mathbf{t}_{d-1,d}^T\mathbf{a}}/2))^{-1} \end{bmatrix}.$$

Next define

$$\hat{\mathbf{a}} := (\boldsymbol{\Gamma}^T\boldsymbol{\Gamma})^{-1}\boldsymbol{\Gamma}^T\hat{\mathbf{q}}. \tag{3.13}$$

Then

$$\sqrt{k}(\hat{\mathbf{a}} - \mathbf{a}) = (\boldsymbol{\Gamma}^T\boldsymbol{\Gamma})^{-1}\boldsymbol{\Gamma}^T\sqrt{k}(\hat{\mathbf{q}} - \boldsymbol{\Gamma}\mathbf{a}) \xrightarrow{d} (\boldsymbol{\Gamma}^T\boldsymbol{\Gamma})^{-1}\boldsymbol{\Gamma}^T\mathbf{b}. \tag{3.14}$$

By solving the equations (3.9) we get

$$\beta_1 := \sqrt{a_1 - \frac{a_2^2}{4a_3}},$$

$$\beta_2 := \sqrt{a_3 - \frac{a_2^2}{4a_1}}$$

and

$$\rho := -\frac{a_2}{2\sqrt{a_1 a_3}}.$$

Cramér's delta method now gives the joint asymptotic normality of the estimators

$$\hat{\beta}_1 := \sqrt{\hat{a}_1 - \frac{\hat{a}_2^2}{4\hat{a}_3}}, \tag{3.15}$$

$$\hat{\beta}_2 := \sqrt{\hat{a}_3 - \frac{\hat{a}_2^2}{4\hat{a}_1}} \tag{3.16}$$

and

$$\hat{\rho} := -\frac{\hat{a}_2}{2\sqrt{\hat{a}_1 \hat{a}_3}}. \tag{3.17}$$

For the proofs of the theorems we need a number of auxiliary results.

LEMMA 3.1. *Suppose for some measure $\nu$ on $\mathbb{R}^d$ and some positive integrable functions $g_1, \ldots, g_d$*

$$\nu([x_1, \infty) \cup \cdots \cup [x_d, \infty)) = \int_{-\infty}^{\infty} \max\{g_1(x - x_1), \ldots, g_d(x - x_d)\}\, dx.$$

*Then*

$$\nu([x_1,\infty) \cap \cdots \cap [x_d,\infty)) = \int_{-\infty}^{\infty} \min\{g_1(x-x_1),\ldots,g_d(x-x_d)\}\,dx.$$

PROOF.    Note that

$$\nu([x_1,\infty) \cup [x_2,\infty) \cup \cdots \cup [x_d,\infty))$$

$$= \sum_{i=1}^{d} \nu([x_i,\infty)) + (-1)\sum_{i \neq j} \nu([x_i,\infty) \cap [x_j,\infty))$$

$$+ \cdots + (-1)^{d-1}\nu([x_1,\infty) \cap [x_2,\infty) \cap \cdots \cap [x_d,\infty)),$$

and for any real $a_1, a_2, \ldots, a_d$,

$$\max(a_1, a_2, \ldots, a_d) = \sum_{i=1}^{d} a_i + (-1)\sum_{i \neq j} a_i \wedge a_j + \cdots + (-1)^{d-1} a_1 \wedge a_2 \wedge \cdots \wedge a_d$$

(both follow by induction). Then the statement also follows by induction.
□

COROLLARY 3.1.    *Hence for our models in* $\mathbb{R}$

$$R_{t_1,t_2,\ldots,t_d}(x_1,x_2,\ldots,x_d) = \int_{-\infty}^{+\infty} \min_i \frac{\phi(s-t_i)}{x_i}\,ds.$$

REMARK 3.2.    We apply Lemma 3.1, for example, to the functions $g_j(x) :=$ $\frac{e^{-\beta|x|}}{w_j}$.

LEMMA 3.2.    *Suppose $p$ is a probability density on $\mathbb{R}$, $p(x) = p(-x)$ for $x > 0$ and $p(x)$ is decreasing for $x > 0$. Then*

(3.18)
$$\int_{-\infty}^{\infty} \min\{p(|s-t_1|),\ldots,p(|s-t_d|)\}\,ds$$
$$= 2\int_{(1/2)(\max_{1 \leq j \leq d} t_j - \min_{1 \leq j \leq d} t_j)}^{\infty} p(s)\,ds.$$

PROOF.    Note that

$$\min\{p(|s-t_1|),\ldots,p(|s-t_d|)\}$$

$$= \min\left\{p\left(\left|s - \min_{1 \leq j \leq d} t_j\right|\right), p\left(\left|s - \max_{1 \leq j \leq d} t_j\right|\right)\right\}$$

$$= \begin{cases} p\left(\left|s - \min_{1 \leq j \leq d} t_j\right|\right), & s > \min_{1 \leq j \leq d} t_j + \frac{1}{2}\left(\max_{1 \leq j \leq d} t_j - \min_{1 \leq j \leq d} t_j\right), \\ p\left(\left|s - \max_{1 \leq j \leq d} t_j\right|\right), & s < \min_{1 \leq j \leq d} t_j + \frac{1}{2}\left(\max_{1 \leq j \leq d} t_j - \min_{1 \leq j \leq d} t_j\right). \end{cases}$$

So the integral on the left-hand side equals

$$2\int_{\min_{1\le j\le d} t_j + (1/2)(\max_{1\le j\le d} t_j - \min_{1\le j\le d} t_j)}^{\infty} p\left(\left|s - \min_{1\le j\le d} t_j\right|\right) ds. \qquad \square$$

COROLLARY 3.2.  *Hence for our models in* $\mathbb{R}$

$$R_{t_1,\ldots,t_d}(1,\ldots,1) = 2\int_{(1/2)(\max_{1\le j\le d} t_j - \min_{1\le j\le d} t_j)}^{\infty} \phi(s)\, ds$$

$$= 2\left\{1 - F\left(\frac{\beta}{2}\left(\max_{1\le j\le d} t_j - \min_{1\le j\le d} t_j\right)\right)\right\}.$$

LEMMA 3.3.  *Under the conditions of Theorem* 3.3,

(3.19)
$$\sqrt{k}\left\{\hat{R}_{t_1,\ldots,t_d}(1,\ldots,1) - 2\left[1 - F\left(\frac{\beta}{2}\left(\max_{1\le j\le d} t_j - \min_{1\le j\le d} t_j\right)\right)\right]\right\}$$
$$\to B_{t_1,\ldots,t_d}(1,\ldots,1)$$

*in distribution with* $B$ *as in Proposition* 3.2.

PROOF.  Combine Proposition 3.2 and Corollary 3.2.  $\square$

LEMMA 3.4.  *Under the conditions of Theorem* 3.4, *for the standard normal and Student models,*

(3.20)  $\sqrt{k}\left\{\hat{R}_{\mathbf{t}_j,\mathbf{t}_m}(1,1) - 2\left(1 - F\left(\frac{\beta}{2}(|\mathbf{t}_j - \mathbf{t}_m|)\right)\right)\right\} \to B_{\mathbf{t}_j,\mathbf{t}_m}(1,1)$

*in distribution. For the exponential model*

$$\sqrt{k}\left\{\hat{R}_{\mathbf{t}_j,\mathbf{t}_m}(1,1) - \frac{1}{2}\left(1 + \frac{\beta}{2}\min(|t_1^{(j)} - t_1^{(m)}|, |t_2^{(j)} - t_2^{(m)}|)\right.\right.$$

(3.21)
$$\left.\left. \times\, e^{-(\beta/2)(|t_1^{(j)} - t_1^{(m)}| + |t_2^{(j)} - t_2^{(m)}|)}\right)\right\}$$

$$\to B_{\mathbf{t}_j,\mathbf{t}_m}(1,1)$$

*in distribution, where* $\mathbf{t}_j = (t_1^{(j)}, t_2^{(j)})$, $\mathbf{t}_m = (t_1^{(m)}, t_2^{(m)})$.

PROOF.  Propositions 2.5, 2.6 and 2.4 give $L_{\mathbf{t}_j,\mathbf{t}_m}(1,1)$, which is $2 - R_{\mathbf{t}_j,\mathbf{t}_m}(1,1)$.  $\square$

PROOF OF THEOREM 3.1.  Follows immediately from statement (3.2) of Proposition 3.2 and Lemma 3.2.  $\square$

PROOF OF THEOREMS 3.3 AND 3.4.  Immediately from Lemmas 3.3, 3.4 and Cramér's delta method.  $\square$

**Acknowledgments.** Comments of two referees made us aware of Smith [14] and helped greatly to improve the presentation.

ECONOMETRIC INSTITUTE
ERASMUS UNIVERSITY
ROTTERDAM
NETHERLANDS
E-MAIL: ldehaan@few.eur.nl

DEPARTMENT OF STATISTICS
AND OPERATIONAL RESEARCH
FACULTY OF SCIENCE
UNIVERSITY OF LISBON
CAMPO GRANDE, BLOCO C6, PISO 4
1749-016 LISBOA
PORTUGAL
E-MAIL: tpereira@fc.ul.pt